\documentstyle[11pt,dgspp]{article}
\newtheorem{theorem}{Theorem}[section]

\newtheorem{proposition}[theorem]{Proposition}
\newtheorem{lemma}[theorem]{Lemma}

\newremark{example}[theorem]{Example}
\newremark{remark}[theorem]{Remark}

\newcommand{\bds}{\bigoplus} 

\newcommand{\dps}{\displaystyle}
\newcommand{\ds}{\oplus} 

\newcommand{\half}{\mbox{$\txs\frac{1}{2}$}}

\newcommand{\Iaut}{I^{a\kern-.05em u\kern-.05em t}}
\newcommand{\Ispl}{I^{s\kern-.055em p\kern-.025em l}}

\newcommand{\lsp}{[\kern-0.15em[} 

\newcommand{\ph}{$p\kern-.1em H\!$}

\newcommand{\rsp}{]\kern-0.15em]} 

\newcommand{\surj}{\rightarrow\kern-.82em\rightarrow}

\newcommand{\txs}{\textstyle}

\newcommand{\II}{\rscr{I{\kern-.55em}I}}

\newcommand{\R}{\Bbb R}
\newcommand{\C}{\Bbb C}

\renewcommand{\l}{\lambda}

\makeatletter
\newcommand{\Ad}[1]{\mathop{\operator@font Ad}\nolimits_{#1}}
\newcommand{\ad}[1]{\mathop{\operator@font ad}\nolimits_{#1}}
\newcommand{\add}[2]{\mathop{\operator@font
   ad}\nolimits^{\dagger}_{#1}{\!#2}}
\newcommand{\Aut}{\mathop{\operator@font Aut}\nolimits}
\newcommand{\diag}{\mathop{\operator@font diag}\nolimits}
\newcommand{\End}{\mathop{\operator@font End}\nolimits}
\renewcommand{\Re}{\mathop{\operator@font Re}\nolimits}
\renewcommand{\Im}{\mathop{\operator@font Im}\nolimits}
\newcommand{\Ric}{\mathop{\operator@font Ric}\nolimits}
\newcommand{\Dspecp}{\mathop{\operator@font {\escr D}spec}\nolimits_\wp}
\newcommand{\specl}{\mathop{\operator@font spec}\nolimits_\ell}
\newcommand{\specp}{\mathop{\operator@font spec}\nolimits_\wp}
\newcommand{\tr}{\mathop{\operator@font tr}\nolimits}
\makeatother

\preprint{JP2}
\title{SKEWADJOINT OPERATORS ON PSEUDOEUCLIDEAN SPACES}
\author{Changrim Jang\thanks{On sabbatical leave; permanent address: 
        School of Mathematics and Applied Physics,
        College of Natural Sciences,
        University of Ulsan,
        Ulsan 680-749,
        Republic of Korea. Email:
        crjang@mail.ulsan.ac.kr}
\quad {\small and} \quad Phillip E. Parker}
\address{Mathematics Department\\
        Wichita State University\\
        Wichita KS 67260-0033\\
        USA\\
        \hspace*{-1em}crjang@math.twsu.edu \qquad\qquad 
phil@math.twsu.edu}
\date{29 January 2003} 
\abstract{
We give a complete classification in canonical forms on finite-dimensional
vector spaces over $\R$.
}
\msc{15A21}{15A57, 53C50}

\begin{document}

\maketitle


\section{Introduction}
\label{intro}
In this paper we give a complete classification of canonical forms for 
skewadjoint operators on pseudoeuclidean vector spaces, complementing 
O'Neill's results \cite{O} on selfadjoint operators.

There is a certain skewadjoint operator in the integration of the geodesic 
equations on a 2-step nilpotent Lie group with a left-invariant metric 
tensor, usually called $J$ after Eberlein \cite{E}. In the Riemannian 
(positive definite) case, $J^2$ is always diagonalizable. This makes 
possible the completely explicit integration of the geodesic equations 
(compare \cite{E} and \cite{CP4}) and the rather complete and explicit 
determination of the conjugate loci (compare \cite{JPk} and \cite{JP1}). 
In attempting to extend some of these results to the pseudoriemannian 
(indefinite) case, it became apparent that we needed a set of canonical 
forms for such skewadjoint operators $J$. Obtaining the list turned out to
be a sufficiently involved process that we are presenting it separately.

Throughout, we consider a (finite-dimensional) vector space $V$ provided 
with an inner product $\langle\,,\rangle$; {\em i.e.,} a nondegenerate, 
symmetric bilinear form which is generally indefinite. When necessary, we 
comment on differences if the form is actually definite. Also, $J:V\to V$ 
will denote a skewadjoint linear operator.

We used \cite{HJ} as our main reference for linear algebra, and \cite{BM} 
to a lesser extent, and refer to them for all details relevant to the 
computational procedures employed.

\section{Preparation for the classification}
\label{prep}

To begin, we consider the minimal polynomial of $J$ in the (preliminary) 
form
$$ p(t) = \prod_{i=1}^p\left( t^2 - 2a_i t + (a_i^2+b_i^2)\right)^{k_i} 
\prod_{i=1}^q\left( t^2 + \l_i^2\right)^{l_i} 
\prod_{i=1}^r\left(t-\mu_i\right)^{m_i} $$
where $a_i\neq 0\neq b_i$, $\l_i >0$, and all $a_i,b_i,\l_i,\mu_i$ are 
real. The types of factors correspond to complex eigenvalues that are not 
pure imaginary, pure imaginary eigenvalues, and real eigenvalues, 
respectively.
\begin{proposition}
For\label{jc1} $\mu_i\neq 0$,
$$\ker\left(J-\mu_i I\right)^{m_i} \mbox{ and \ } \ker\left(J^2 - 2a_i J +
(a_i^2+b_i^2)I\right)^{k_i}$$
are null subspaces of $V$.
\end{proposition}
\begin{proof}
If $u,v \in \ker(J-\mu_i I)^{m_i}$, then
\begin{eqnarray*}
\langle u,(J+\mu_i I)^{m_i}v\rangle &=& \langle (-J+\mu_i 
   I)^{m_i}u,v\rangle\\
&=& \langle (-1)^{m_i}(J-\mu_i I)^{m_i}u,v\rangle\\
&=& \langle 0,v\rangle\\
&=& 0\,.
\end{eqnarray*}
Note that $(J+\mu_i I)^{m_i}:\ker(J-\mu_i I)^{m_i} \to \ker(J-\mu_i
I)^{m_i}$ is nonsingular.  Thus it follows that $\langle u,v\rangle=0$ for
all $u,v$ and $\ker(J-\mu_i I)^{m_i}$ is null as claimed.  The other case
is done similarly.
\end{proof}
\begin{proposition}
For\label{jc2} nonzero $\mu_i$, $\ker(J+\mu_i I)^{m_i}$ and $\ker(J^2 + 2a_i 
J + (a_i^2+b_i^2)I)^{k_i}$ are complementary null subspaces to 
$\ker(J-\mu_i I)^{m_i}$ and to $\ker(J^2 - 2a_i J + (a_i^2+b_i^2)I)$, 
respectively.
\end{proposition}
\begin{proof}
Suppose $\ker(J+\mu_i I) = 0$. Then we see that $J+\mu_i I: W\to W$ is 
nonsingular, where
\begin{eqnarray*}
W &=& \bds_{i=1}^p \ker\left(J^2 -2a_i J+(a_i^2+b_i^2)I\right)^{k_i} \\
&& {}+\, \bds_{i=1}^q \ker\left(J^2+\l_i I\right)^{l_i} +\, \bds_{j\neq 
   i}\ker\left( J-\mu_j I\right)^{m_j}
\end{eqnarray*}
By the previous result, if $0\neq u\in\ker(J-\mu_i I)$, then there exists 
$w\in W$ such that $\langle u,w\rangle =1$. Since $J+\mu_i I$ is 
nonsingular on $W$, there exists $w'\in W$ such that $(J+\mu_i I)w'=w$. 
But now we have
\begin{eqnarray*}
1\, =\, \langle u,w\rangle &=& \langle u,(J+\mu_i I)w'\rangle\\
&=& \langle(-J+\mu_i I)u,w'\rangle\\
&=& -\langle(J-\mu_i I)u,w'\rangle\\
&=& 0
\end{eqnarray*}
which is a contradiction. Thus a power term $(t+\mu_i)^{m'_i}$ must appear
in $p(t)$. We shall show that $m'_i=m_i$.

Suppose $m'_i<m_i$. We can choose $u\in V$ such that $(J-\mu_i I)^{m_i-1}u
\neq 0$ and $(J-\mu_i I)^{m_i}u=0$. Then there exists $v\in W$ such that 
$\langle (J-\mu_i I)^{m_i-1}u,v\rangle=1$. Since $m_i-1\ge m'_i$, 
$(-1)^{m_i-1}(J+\mu_i I)^{m_i-1}v\in W_1$ where
\begin{eqnarray*}
W_1 &=& \bds_{i=1}^p \ker\left(J^2-2a_i J+(a_i^2+b_i^2)I \right)^{k_i}\\
&& {}+\, \bds_{i=1}^q \ker\left(J^2+\l_i I \right)^{l_i} + 
   \!\bds_{\mu_j\neq\pm\mu_i}\!\ker\left( J-\mu_j I\right)^{m_j}
\end{eqnarray*}
Since $(J+\mu_i I)^{m_i} : W_1\to W_1$ is nonsingular, there exists $v'\in
W_1$ such that $(J+\mu_i)^{m_i}v' = (-1)^{m_i-1}(J+\mu_i I)^{m_i-1}v$. 
But now we have
\begin{eqnarray*}
1\, =\, \langle u,(-1)^{m_i-1}(J+\mu_i I)^{m_i-1}v\rangle &=& \langle 
u,(J+\mu_i I)^{m_i}v'\rangle\\
&=& \langle(-1)^{m_i}(J-\mu_i I)^{m_i}u,v'\rangle\\
&=& 0
\end{eqnarray*}
which is a contradiction. Thus $m'_i\ge m_i$. Interchanging roles, we get 
$m_i\ge m'_i$ so $m'_i=m_i$.

Now we show that $\ker(J+\mu_i I)^{m_i}$ is a complementary null subspace 
to $\ker(J-\mu_i I)^{m_i}$. For $0\neq u\in \ker(J-\mu_i I)^{m_i}$, there 
exists $w\in W$ such that $\langle u,w\rangle=1$. Now write $w=w_1+w_2$ 
where $w_1\in \ker(J+\mu_i I)^{m_i}$ and $w_2\in W_1$. Since $(J+\mu_i 
I)^{m_i} : W_1\to W_1$ is nonsingular, there exists $w'_2$ such that 
$(J+\mu_i I)^{m_i}w'_2 = w_2$. So we have
\begin{eqnarray*}
1\, =\, \langle u,w\rangle = \langle u,w_1+w_2\rangle &=& \langle u,w_1 
\rangle + \langle u,w_2\rangle \\
&=& \langle u,w_1\rangle + \langle u,(J+\mu_i I)^{m_i}w'_2\rangle\\
&=& \langle u,w_1\rangle + \langle (-1)^{m_i}(J-\mu_i I)^{m_i}u,w'_2 
   \rangle\\
&=& \langle u,w_1\rangle + 0\\
&=& \langle u,w_1\rangle\,.
\end{eqnarray*}
It follows that $\ker(J-\mu_i I)^{m_i}$ and $\ker(J+\mu_i I)^{m_i}$ are 
complementary null subspaces (of each other).

Similarly, one may show that $\ker(J^2 + 2a_i J + (a_i^2+b_i^2)I)$ and 
$\ker(J^2 - 2a_i J + (a_i^2+b_i^2)I)$ are likewise complementary.
\end{proof}

{From} these two propositions it now follows that the minimal polynomial 
of $J$ can be written in the form
\begin{equation}
\label{mp} \begin{array}{rcl}
\dps p(t) &=& \dps\prod_{i=1}^p \left(t^2 - 2a_i t + (a_i^2+b_i^2)
   \right)^{k_i} \left(t^2 + 2a_i t + (a_i^2+b_i^2) \right)^{k_i}\\[2.5ex]
&& {} \dps\cdot \prod_{i=1}^q \left(t^2+\l_i^2\right)^{l_i} \cdot 
   \prod_{i=1}^r \left(t-\mu_i\right)^{m_i}\left(t+\mu_i\right)^{m_i} 
   \cdot t^s \end{array}
\end{equation}
where  $a_i\neq 0\neq b_i$, $\l_i >0$, $\mu_i\neq 0$, and all
$a_i,b_i,\l_i,\mu_i$ are real.

\begin{proposition}
\label{jc3}$\ker(J^2 - 2a_i J + (a_i^2+b_i^2)I)^{k_i}\ds \ker(J^2 + 2a_i J + 
(a_i^2+b_i^2)I)^{k_i}$ for $1\le i\le p$, $\ker(J^2+\l_i^2 I)^{l_i}$ for 
$1\le i\le q$, $\ker(J-\mu_i I)^{m_i}\ds\ker(J+\mu_i I)^{m_i}$ for $1\le 
i\le r$, and $\ker J^s$ are $J$-invariant, mutually orthogonal,
nondegenerate subspaces of $V$.
\end{proposition}
\begin{proof}
These are rather straightforward calculations.  As an example, we show
that each $\ker(J^2+\l_i^2)^{l_i}$ is nondegenerate. 

Suppose a vector $u\in \ker(J^2+\l_i^2 I)^{l_i}$ satisfies $\langle 
u,v\rangle=0$ for all $v\in \ker(J^2+\l_i^2 I)^{l_i}$. Let
\begin{eqnarray*}
W &=& \bds_{i=1}^p \left( \ker(J^2-2a_i J+(a_i^2+b_i^2)I)^{k_i} \ds 
   \ker(J^2+2a_i J+(a_i^2+b_i^2)I)^{k_i} \right) \\
&& \ds\ \bds_{i=1}^r \left( \ker(J-\mu_i I)^{m_i} \ds \ker(J+\mu_i 
I)^{m_i}
   \right) \\
&& \ds\ \bds_{j\neq i}\ker(J^2+\l_j^2 I)^{l_j} \ds \ker J^s\,.
\end{eqnarray*}
Then $(J^2+\l_i^2 I)^{l_i} :W\to W$ is nonsingular, hence surjective. This
and $\langle u,(J^2+\l_i^2 I)^{l_i}w\rangle = \langle (J^2+\l_i^2 
I)^{l_i}u,w\rangle = 0$ for all $w\in W$ imply that $\langle u,w\rangle = 
0$ for all $w\in W$, thus $\langle u,v\rangle=0$ for all $v\in V$, so 
$u=0$ since $V$ is nondegenerate. Consequently $\ker(J^2+\l_i^2 I)^{l_i}$ 
is also nondegenerate.
\end{proof}

It now follows that $J$ can be decomposed as an orthogonal direct sum of 
skewadjoint operators on each of these subspaces:
\begin{eqnarray*}
&& \ker(J^2-2a_i J+(a_i^2+b_i^2)I)^{k_i} \ds \ker(J^2+2a_i 
   J+(a_i^2+b_i^2)I)^{k_i}\\
&& \hspace{25em}\mbox{for }1\le i\le p\, ;\\
&& \ker(J^2+\l_i^2 I)^{l_i} \quad\mbox{for } 1\le i\le q\, ;\\
&& \ker(J-\mu_i I)^{m_i} \ds \ker(J+\mu_i I)^{m_i} \quad\mbox{for } 1\le 
   i\le r\, ;\\
&& \ker J^s .
\end{eqnarray*}
Thus it will suffice to classify skewadjoint operators of each type.

\section{Minimal polynomial \boldmath$(t^2-2at+(a^2+b^2))^k(t^2+2at+(a^2+b^2))^k$}
\label{ce}

If $J$ is such a skewadjoint operator on $V$, then we can also consider 
$J$ as a skewadjoint operator on the complexification $U = V^{\C}$ of $V$.
Here the inner product is extended by bilinearity in the usual way.
$$\langle v_1+iw_1, v_2+iw_2\rangle = \langle v_1,v_2\rangle - \langle 
w_1,w_2\rangle + i(\langle w_1,v_2\rangle + \langle v_1,w_2\rangle ) $$

The minimal polynomial of $J$ factors over $\C$ as
$$ (t-(a+ib))^k (t-(a-ib))^k (t+(a+ib))^k (t+(a-ib))^k\, .$$
Then there exists a vector $u+iv\in U$ such that $(J-(a+ib)I)^{k-1}(u+iv) 
\neq 0$ and $(J-(a+ib)I)^k(u+iv) = 0$. Since $\Re (J-(a+ib)I)^{k-1}(u+iv)$
and $\Im (J-(a+ib)I)^{k-1}(u+iv)$ lie in $\ker (J^2-2aJ+(a^2+b^2)I)^k$, 
there exist $w,x\in \ker(J^2+2aJ+(a^2+b^2)I)^k$ such that
\begin{equation}
\label{j3.1}
\langle (J-(a+ib)I)^{k-1}(u+iv), w+ix\rangle = 2
\end{equation}
by Proposition \ref{jc2}. Now we need the following result.
\begin{lemma}
If\label{lem'} $u+iv\in\ker(J-(a+ib)I)^k$ and $w+ix\in\ker(J+(a-ib)I)^k$, 
then $\langle u+iv,w+ix\rangle = 0$.
\end{lemma}
\begin{proof}
Since 
$$\left( J+(a+ib)I\right)^k : \ker\left( J+(a-ib)I\right)^k \to 
   \ker\left( J+(a-ib)I\right)^k$$ 
is nonsingular, there exists $w'+ix'\in\ker\left( J+(a-ib)I\right)^k$
such that 
$$\left( J+(a+ib)I\right)^k (w'+ix') = w+ix\,.$$
Thus,
\begin{eqnarray*}
\langle u+iv,w+ix\rangle &=& \langle u+iv, \left(J+(a+ib)I \right)^k
   (w'+ix') \rangle \\
&=& \langle (-1)^k\left( J-(a+ib)I\right)^k(u+iv), w'+ix'\rangle \\
&=& \langle 0,w'+ix'\rangle\\
&=& 0\,.\eop
\end{eqnarray*}
\end{proof}
Write $w+ix = w_1 + ix_1 + w_2 + ix_2$ with $w_1+ix_1 \in \ker\left( 
J+(a+ib)I\right)^k$ and $w_2+ix_2 \in \ker\left( J+(a-ib)I\right)^k$. Then
from (\ref{j3.1}) and Lemma \ref{lem'} we obtain
\begin{eqnarray*}
2 &=& \langle (J-(a+ib)I)^{k-1}(u+iv),w+ix\rangle\\
&=& \langle (-(a+ib)I)^{k-1}(u+iv),w_1+ix_1+w_2+ix_2\rangle\\
&=& \langle (J-(a+ib)I)^{k-1}(u+iv),w_1+ix_1\rangle
\end{eqnarray*}
This implies that we may assume that 
$$w+ix \in \ker\left( J+(a+ib)I\right)^k$$ 
in (\ref{j3.1}). Since $w-ix\in\ker(J+(a-ib)I)^k$, we have {\em via\/} 
Lemma \ref{lem'}
\begin{equation}
\label{j3.2}
\langle (J-(a+ib)I)^h(u+iv), (J+(a-ib)I)^j(w-ix)\rangle = 0
\end{equation}
for all nonnegative integers $h,j$.

Also, we can find complex numbers $a_1,a_2,\ldots,a_{k-1}$ such that
\begin{eqnarray*}
\langle (u+iv) + a_1(J-(a+ib)I)(u+iv) + \cdots \hspace{5em} &&\\
+ a_{k-1}(J-(a+ib)I)^{k-1}(u+iv), w+ix \rangle &=& 0\, ,\\
\langle (J-(a+ib)I)(u+iv) + a_1(J-(a+ib)I)^2(u+iv) + \cdots &&\\
+ a_{k-2}(J-(a+ib)I)^{k-1}(u+iv), w+ix \rangle &=& 0\, ,\\
\vdots\hspace{10em}\vdots\hspace{10em} &\vdots & \\
\langle (J-(a+ib)I)^{k-2}(u+iv) + a_1(J-(a+ib)I)^{k-1}(u+iv), w+ix \rangle
&=& 0\,.
\end{eqnarray*}
Replacing $u+iv$ by 
$$(u+iv) + a_1(J-(a+ib)I)(u+iv) + \cdots + 
   a_{k-1}(J-(a+ib)I)^{k-1}(u+iv)$$
we have
\begin{equation}
\label{j3.3}
\langle (J-(a+ib)I)^h(u+iv), w+ix\rangle = 0 \,\mbox{ for }\, 0\le h\le 
   k-2\,.
\end{equation}
Note that (\ref{j3.1}) and (\ref{j3.2}) continue to hold as well.

Let $u_1 = u$ and $v_1=v$ and set
\begin{equation}
\label{j3.4}
(J-(a+ib)I)^h(u_1+iv_1) = u_{h+1}+iv_{h+1}\,\mbox{ for }\, 1\le h\le 
   k-1\,.
\end{equation}
Then
\begin{equation}
\label{j3.5}
\begin{array}{rcl} Ju_h &=& au_h-bv_h+u_{h+1}\\[.5ex]
Jv_h &=&bu_h+av_h+v_{h+1} \end{array}\quad\mbox{for }\, 1\le h\le k-1
\end{equation}
\begin{equation}
\label{j3.6}
\begin{array}{rcl} Ju_k &=& au_k-bv_k\\[.5ex]
Jv_k &=& bu_k+av_k \end{array}\hspace{11.3em}
\end{equation}
Now (\ref{j3.1}) implies that 
\begin{eqnarray*}
\left(J+(a+ib)I\right)^{k-1}(w+ix) &\neq& 0\, ,\\ 
\left(J+(a+ib)I\right)^{k}(w+ix) &=& 0\,. 
\end{eqnarray*}
Let $w_1=w$ and $x_1=x$ and set
\begin{equation}
\label{j3.7}
(J+(a+ib)I)^h(w_1+ix_1) = w_{h+1}+ix_{h+1}\,\mbox{ for }\, 1\le h\le 
   k-1\,.
\end{equation}
Then
\begin{equation}
\label{j3.8}
\begin{array}{rcl} Jw_h &=& -aw_h+bx_h+w_{h+1}\\[.5ex]
Jx_h &=& -bw_h-ax_h+x_{h+1} \end{array}\quad\mbox{for }\, 1\le h\le k-1
\end{equation}
\begin{equation}
\label{j3.9}
\begin{array}{rcl} Jw_k &=& -aw_k+bx_k\\[.5ex]
Jx_k &=& -bw_k-ax_k \end{array}\hspace{11.4em}
\end{equation}

By (\ref{j3.1}), (\ref{j3.2}), and (\ref{j3.3}),
$$ \langle (J-(a+ib)I)^h(u_1+iv_1), (-1)^j(J+(a+ib)I)^j(w_1+ix_1) 
   \rangle $$
$$ {} = \left\{ \begin{array}{rl} 2 & \mbox{ if }\, h+j=k-1\, ,\\
   0 & \mbox{ otherwise,} \end{array} \right. $$
$$ \langle (J-(a+ib)I)^h(u_1+iv_1), (-1)^j(J+(a-ib)I)^j(w_1-ix_1) 
   \rangle =0 $$
for all nonnegative integers $h,j$. This implies that
$$ \begin{array}{rcl}
   \langle u_{h+1}+iv_{j+1}, (-1)^j(w_{j+1}+ix_{j+1}) \rangle &=& \left\{ 
   \begin{array}{rl} 2 & \mbox{ if }\, h+j=k-1\, ,\\
   0 & \mbox{ otherwise,} \end{array} \right. \\[3ex]
   \langle u_{h+1}+iv_{h+1},w_{j+1}-ix_{j+1} \rangle &=& 0 \end{array} $$
for all nonnegative integers $h,j$. Thus we have all inner products zero
except
$$ \begin{array}{rcl}
   \langle u_h,w_j\rangle &=& (-1)^{j-1}\\[.5ex]
   \langle v_h,x_j \rangle &=& (-1)^j \end{array}\quad\mbox{ if }\, 
   h+j=k+1\,.$$
(Recall that $\ker(J^2-2aJ+(a^2+b^2)I)$ and $\ker(J^2+2aJ+(a^2+b^2)I)$ are
null subspaces.)
Note that $\{ u_1, v_i, \ldots , u_k, v_k, w_1, x_1, \ldots , w_k, x_k \}$
spans a $J$-invariant, nondegenerate subspace. Hence there exists a 
$J$-invariant, nondegenerate subspace $N$ of $V$ such that 
$$\lsp u_1, v_1, \ldots , u_k, v_k, w_1, x_1, \ldots , w_k, x_k \rsp\ds 
   N$$
is an orthogonal direct sum.  On the subspace $\lsp u_1, v_1, \ldots
, w_k, x_k \rsp$, the matrix of $J$ is $2k\times 2k$ of the following 
form.
$$ \arraycolsep .325em \left[ 
\begin{array}{cccccccccccccccccc} a & b &&&&&&&&&&&&&&&&\\ 
-b & a &&&&&&&&&&&&&&&&\\ 1 & 0 & a & b &&&&&&&&&&&&&&\\
0 & 1 & -b & a &&&&&&&&&&&&&&\\[-.6ex] &&&& \ddots &&&&&&&&&&&&&\\
&&&&& a & b &&&&&&&&&&&\\ &&&&& -b & a &&&&&&&&&&&\\
&&&&& 1 & 0 & a & b &&&&&&&&&\\ &&&&& 0 & 1 & -b & a &&&&&&&&&\\
&&&&&&&&& -a & -b &&&&&&&\\ 
&&&&&&&&& b & -a &&&&&&&\\ &&&&&&&&& 1 & 0 & -a & -b &&&&&\\
&&&&&&&&& 0 & 1 & b & -a &&&&&\\[-.6ex] &&&&&&&&&&&&& \ddots &&&&\\
&&&&&&&&&&&&&& -a & -b &&\\ &&&&&&&&&&&&&& b & -a &&\\
&&&&&&&&&&&&&& 1 & 0 & -a & -b \\ &&&&&&&&&&&&&& 0 & 1 & b & -a \\
\end{array} \right] $$
Continuing this process, we get an orthogonal direct sum of $J$-invariant
subspaces $V_h$ such that on each of them, $J$ has a matrix of the
preceding form but of size $2r\times 2r$ for some $r$ with $1\le r\le k$,
with respect to a basis of the specified type.

\section{Minimal polynomial \boldmath$(t^2+\l^2)^l$}
\label{pie}

Let $J$ be a skewadjoint operator on $V$ with such a minimal polynomial.
\begin{lemma}
\label{lem}
There exists a vector $v \in V$ such that $v$, $(J^2+\l^2 I)v$, \ldots , 
$(J^2+\l^2 I)^{l-1}v$ are linearly independent and $\langle v, (J^2+\l^2 
I)^{l-1}v \rangle \neq 0$.
\end{lemma}
\begin{proof}
Since the minimal polynomial of $J$ is $(t^2+\l^2)^l$, there exists $v\in 
V$ such that $v$, $(J^2+\l^2 I)v$, \ldots , $(J^2+\l^2 I)^{l-1}v$ are 
linearly independent.

If $\langle v, (J^2+\l^2 I)^{l-1}v \rangle = 0$, then there exists $w\in 
V$ such that $\langle w, (J^2+\l^2 I)^{l-1}v \rangle \neq 0$. If $\langle 
w, (J^2+\l^2 I)^{l-1}w \rangle \neq 0$, then we simply replace $v$ with 
$w$ and verify that we have achieved the desired result; otherwise, 
replace $v$ with $v+w$ and verify that this works.
\end{proof}

So we may start with $v\in V$ as in the Lemma. Set $(J+i\l I)^l v = 
u_1+iv_1$ and observe that $(u_1+iv_1)$, $(J-i\l I)(u_1+iv_1)$, \ldots , 
$(J-i\l I)^{l-1}(u_1+iv_1)$ are linearly independent over $\C$ and that 
$(J-i\l I)^l (u_1+iv_1) = (J^2+\l^2 I)^l v = 0$. Let
$$ \left( J-i\l I\right)^h (u_1+iv_1) = u_{h+1}+iv_{h+1}\,\mbox{ for }\, 
   0\le h\le l-1\,.$$
First off, we have
$$ \langle u_h+iv_h, u_j+iv_j\rangle = 0\,\mbox{ for }\, 1\le h,j\le l$$
since
\begin{eqnarray*}
\lefteqn{\langle u_h+iv_h, u_j+iv_j\rangle }\hspace{3.4em}\\
&=& \langle (J-i\l I)^{h-1}(J+i\l I)^l v, (J-i\l I)^{j-1}(J+i\l I)^l v 
   \rangle \\
&=& (-1)^l\langle (J-i\l I)^{h-1}(J-i\l I)^l(J+i\l I)^l v, (J-i\l 
   I)^{j-1}v \rangle \\
&=& (-1)^l\langle (J-i\l I)^{h-1}(J^2+\l^2 I)^l v, (J-i\l I)^{j-1}v 
   \rangle \\
&=& 0\,.
\end{eqnarray*}
This implies that
\begin{equation}
\label{c1}
\begin{array}{rcl}
\dps\langle u_h, u_j \rangle - \langle v_h, v_j \rangle &=& 0\\[1ex]
\dps\langle u_h, v_j \rangle + \langle v_h, u_j \rangle &=& 0 
\end{array}
\hspace{2em}\mbox{for }\, 1\le h,j\le l\,.
\end{equation}

Next, we have
\begin{eqnarray*}
\lefteqn{\langle u_l+iv_l, u_h-iv_h \rangle }\hspace{3em}\\
&=& \langle (J-i\l I)^{l-1}(J+i\l I)^l v, (J+i\l I)^{h-1}(J-i\l I)^l v 
   \rangle \\
&=& (-1)^{h-1}\langle (J-i\l I)^{l+h-2}(J+i\l I)^l v, (J-i\l I)^l v 
   \rangle \\
&=& (-1)^{h-1}\langle (J-i\l I)^{h-2}(J^2+\l^2 I)^l v, (J-i\l I)^l v 
   \rangle
\end{eqnarray*}
whence
\begin{equation}
\label{c2}
\langle u_l+iv_l, u_h-iv_h \rangle = 0\hspace{2em}\mbox{for }\, 2\le h\le 
   l\,.
\end{equation}

Note that $J^2(J^2+i\l I)^{l-1}v = -\l^2(J^2+i\l^2 I)^{l-1}v$. We use this
in the following computation.
\begin{lemma}
If\label{lem''} $l$ is odd, then $\langle u_1+iv_1, u_l-iv_l \rangle$ is a 
nonzero real number, and if $l$ is even then it is nonzero pure imaginary.
Moreover, $\langle u_h+iv_h, u_j-iv_j \rangle $ is real when $h+j$ is even
and pure imaginary when $h+j$ is odd.
\end{lemma}
\begin{proof} We compute
\begin{eqnarray*}
\lefteqn{ \langle u_1+iv_1, u_l-iv_l \rangle }\\
&=& \langle (J+i\l I)^l v, (J+i\l I)^{l-1}(J-i\l I)^l v \rangle \\
&=& -\langle (J+i\l I)^{l+1}v, (J^2+\l^2 I)^{l-1}v \rangle \\
&=& -\langle \left({\txs J^{l+1} + {l+1\choose 1}(i\l)J^l + \cdots + 
   {l+1\choose l}(i\l)^l J + (i\l)^{l+1}I }\right) v, (J^2+\l^2 I)^{l-1}v 
   \rangle \,.
\end{eqnarray*}
If $l$ is odd,
\begin{eqnarray*}
\langle u_1+iv_1, u_l-iv_l \rangle 
&=& -\langle \left({\txs J^{l+1} - {l+1\choose 2}\l^2 J^{l-1} + 
   {l+1\choose 4}\l^4 J^{l-3} - \cdots }\right. \\
&& \hspace{5em}\left.{\txs {} + (-1)^{\frac{l+1}{2}} \l^{l+1}I}\right) v, 
   (J^2+\l^2 I)^{l-1} \rangle \\
&=& -\langle \left( {\txs (-1)^{\frac{l+1}{2}}\l^{l+1} - {l+1\choose 2} 
   \l^2(-1)^{\frac{l-1}{2}}\l^{l-1} + \cdots }\right. \\
&& \hspace{5em}\left.{\txs {} + (-1)^{\frac{l+1}{2}}\l^{l+1}}\right) v, 
   (J^2+\l^2 I)^{l-1}v \rangle \\
&=& (-1)^{\frac{l+3}{2}}\l^{l+1} \left( {\txs 1 + {l+1\choose 2} + 
   {l+1\choose 4} + \cdots }\right. \\
&& \hspace{5em}\left.{\txs {} + {l+1\choose l-1} + 1 }\right) \langle 
   v, (J^2+\l^2 I)^{l-1}v \rangle
\end{eqnarray*}
which is a nonzero real number. If $l$ is even, an analogous computation 
yields a nonzero pure imaginary number. The rest is done similarly.
\end{proof}

Now assume that $l$ is odd and consider a vector of the form
$$ w+ix = (u_1+iv_1) + \alpha_1(u_2+iv_2) + \cdots + 
   \alpha_{l-1}(u_l+iv_l)\,.$$
Then we can determine the (complex) coefficients $\alpha_i$ such that
$$\langle w+ix, (J+i\l I)^h(w-ix)\rangle = 0 \hspace{2em}\mbox{for }\, 
0\le h\le l-2\,.$$
Indeed, consider
\begin{eqnarray*}
\lefteqn{ \langle w+ix, (J+i\l I)^{l-2}(w-ix)\rangle }\hspace{3em}\\
&=& \langle (u_1+iv_1) + \alpha_1(u_2+iv_2) + \cdots + 
   \alpha_{l-1}(u_l+iv_l), \\
&& \hspace{3em}(u_{l-1}-iv_{l-1}) + \bar\alpha_1(u_l-iv_l) \rangle \\
&=& \langle u_1+iv_1, u_{l-1}-iv_{l-1} \rangle + \alpha_1\langle u_2+iv_2,
   u_{l-1}-iv_{l-1} \rangle \\
&& \hspace{3em}{} + \bar\alpha_1\langle u_1+iv_1, u_l-iv_l 
   \rangle\,.
\end{eqnarray*}
Since $\langle u_1+iv_1, u_{l-1}-iv_{l-1}\rangle$ is pure imaginary and 
$\langle u_2+iv_2, u_{l-1}-iv_{l-1}\rangle = -\langle u_1+iv_1, u_l-iv_l 
\rangle$ is real, we can find a pure imaginary $\alpha_1$ such that
$$\langle w+ix, (J+i\l I)^{l-2}(w-ix)\rangle = 0\,.$$
At the next step, there are six such inner product terms (with 
coefficients) which are all real. We obtain a real $\alpha_2$ such that
$$\langle w+ix, (J+i\l I)^{l-3}(w-ix) \rangle = 0\,.$$
Continuing inductively, we find all the desired coefficients.

Replacing the previous $u_1+iv_1$ with $w+ix$ so determined, using 
(\ref{c1}), (\ref{c2}), and Lemma \ref{lem''}, and rescaling, we have
$$ \begin{array}{rcl}
\dps\langle u_h+iv_h, u_j+iv_j \rangle &=& 0 \qquad\mbox{for }\, 1\le 
   h,j\le l\, ;\\[1ex]
\dps\langle u_h+iv_h, u_j-iv_j \rangle &=& \left\{ 
   \begin{array}{cl}
   (-1)^j(\pm 2) & \quad\mbox{if }\, h+j=l+1\, ,\\[.5ex]
   0 & \quad\mbox{otherwise.} \end{array} \right. \end{array} $$
In other words, all inner products among the $u_h$ and $v_j$ are zero 
except
$$\langle u_h, u_j \rangle = \langle v_h, v_j \rangle = (-1)^j(\pm 1) 
   \quad\mbox{if }\, h+j=l+1\,.$$
Since $(J-i\l I)(u_h+iv_h) = u_{h+1}+iv_{h+1}$ for $1\le h\le l-1$ and 
$(J-i\l I)(u_l+iv_l) = 0$, the matrix of $J$ on $\lsp u_1,v_1,\ldots 
,u_l,v_l \rsp$ is $2l\times 2l$ of this form.
$$ \arraycolsep .4em \left[ 
\begin{array}{ccccccccc} 0 & \l &&&&&&&\\ 
-\l & 0 &&&&&&&\\ 1 & 0 & 0 & \l &&&&&\\
0 & 1 & -\l & 0 &&&&&\\[-.6ex] &&&& \ddots &&&&\\
&&&&& 0 & \l &&\\ &&&&& -\l & 0 &&\\
&&&&& 1 & 0 & 0 & \l \\ &&&&& 0 & 1 & -\l & 0 \\
\end{array} \right] $$

If $l$ is even, the condition on inner products of the basis vectors 
becomes
$$\langle v_h,u_j \rangle = -\langle u_h,v_j \rangle = (-1)^j(\pm 1) 
\quad\mbox{if }\, h+j=l+1$$
with all others vanishing as before, while the $2l\times 2l$ form of the 
matrix of $J$ remains the same.

We finish by continuing this process as at the end of the preceding 
section, with $1\le r\le l$ now.

\section{Minimal polynomial \boldmath$(t-\mu)^m(t+\mu)^m$}
\label{re}

Let $J:V\to V$ be a skewadjoint operator with such a minimal polynomial. 
Then $V=\ker(J-\mu I)^m \ds \ker(J+\mu I)^m = W_1\ds W_2$ where $W_1$ and 
$W_2$ are complementary null subspaces by Propositions \ref{jc1} 
and \ref{jc2}. There exists $v\in W_1$ such that $v, (J-\mu I)v, \ldots , 
(J-\mu I)^{m-1}v$ are linearly independent and $(J-\mu I)^m v=0$. Then we 
can choose $w\in W_2$ such that $\langle(J-\mu I)^{m-1}v,w\rangle=1$ and 
there exist real numbers $a_1,\ldots,a_{m-1}$ such that
\begin{eqnarray*}
\langle v+a_1(J-\mu I)v + \cdots +a_{m-1}(J-\mu I)^{m-1}v,w\rangle &=& 0\\
\langle (J-\mu I)v + a_1(J-\mu I)^2 v +\cdots+a_{m-2}(J-\mu I)^{m-1} v,w
   \rangle &=& 0\\
\vdots\hspace{10em}\vdots\hspace{5em} && \,\vdots\\
\langle (J-\mu I)^{m-2}v + a_1(J-\mu I)^{m-1}v,w\rangle &=& 0
\end{eqnarray*}

If we replace $v$ by $v+a_1(J-\mu I)v+\cdots+a_{m-1}(J-\mu I)^{m-1}v$ and 
set $v_1=v$, $v_2=(J-\mu I)v_1$, \ldots , $v_m = (J-\mu I)v_{m-1}$, and 
$w_1=w$, $w_2=(J+\mu I)w_1$, \ldots , $w_m=(J+\mu I)w_{m-1}$, then all 
inner products are zero except
\begin{eqnarray*}
\langle v_i,w_j\rangle &=& \langle(J-\mu I)^{i-1}v,(J+\mu 
   I)^{j-1}w\rangle\\
&=& \langle(-1)^{j-1}(J-\mu I)^{i+j-2}v,w\rangle\\
&=& (-1)^{j-1}\quad\mbox{if }\, i+j=m+1.
\end{eqnarray*}

We also have a $J$-invariant, nondegenerate subspace $N$ such that $V = 
\lsp v_1,\ldots ,v_m, w_1,\ldots ,w_m \rsp \ds N$ is an orthogonal direct 
sum. With respect to the obvious basis, the matrix of $J$ on the first 
subspace is block-diagonal
$$ \arraycolsep .4em
\left[ \begin{array}{cccccccc} \mu &&&&&&&\\ 1 & \mu &&&&&&\\ & 1 & 
\ddots &&&&&\\ && 1 & \mu &&&&\\ &&&& -\mu &&&\\ &&&& 1 & -\mu &&\\ &&&&& 
1 & \ddots &\\ &&&&&& 1 & -\mu \end{array} \right] $$
with two elementary Jordan diagonal blocks, each $m\times m$.

We finish as before, with $1\le r\le m$ now.

\section{Minimal polynomial \boldmath$t^s$}
\label{nil}

Let $J$ be such a skewadjoint operator. Then there exists a subspace $V_1$
of $V$ such that $V =\ker J^{s-1}\ds V_1$, whence $J^{s-1}V = J^{s-1}V_1$.
Since $J^{s-1}V_1$ is a null subspace, there exists a complementary null 
subspace $W_1$ and a subspace $N$ which is orthogonal to their direct sum 
so that $V = \left( J^{s-1}V_1 \ds W_1 \right) \ds N$. Then $J^{s-1}V = 
J^{s-1}V_1 = J^{s-1}W_1 \ds J^{s-1}N$. Since $W_1$ is a complementary null
subspace to $J^{s-1}V_1$, this implies that $J^{s-1}V_1 = J^{s-1}W_1$ and 
that $J^{s-1}N \subseteq J^{s-1}W_1$. Thus we may assume that $W_1 = V_1$
and we have $$V = J^{s-1}V_1 \ds V_1 \ds N$$ where $J^{s-1}V_1$ and $V_1$ 
are complementary null subspaces (of each other) and $N$ is orthogonal to 
their (direct) sum.

Firstly, assume that $s$ is even. Consider $J^{s-1}v_1$ for a nonzero
vector in $V_1$. Then there exists $w_1\in V_1$ such that $\langle 
J^{s-1}v_1,w_1 \rangle = 1$ and $v_1$ and $w_1$ are linearly independent.
\begin{lemma}
$\lsp v_1, \ldots ,J^{s-1}v_1, w_1, \ldots ,J^{s-1}w_1 \rsp$ is a 
$J$-invariant, nondegenerate, $2s$-dimensional subspace of $V$.
\end{lemma}
\begin{proof}
Suppose $a_1v_1 + \cdots + a_s J^{s-1}v_1 + b_1w_1 + \cdots + 
b_sJ^{s-1}w_1 = 0$ for real $a_i,b_i$. Applying $J^{s-1}$ to this 
equation, we obtain $a_1J^{s-1}v_1 + b_1J^{s-1}w_1 = 0$. Then $\langle 
a_1J^{s-1}v_1 + b_1J^{s-1}w_1, w_1 \rangle = a_1 = 0$ whence $b_1=0$. 
Continuing this process, all the coefficients vanish so our list is a 
basis and the subspace is $2s$-dimensional.

If $\langle a_1v_1 + \cdots + J^{s-1}v_1 + b_1w_1 + \cdots + J^{s-1}w_1, u
\rangle = 0$ for all $u$ in our subspace with real $a_i,b_i$, then
replacing $u$ by $J^{s-1}w_1$ yields $a_1=0$.  Continuing in this way, all
the coefficients vanish again; this time it shows that our subspace is
nondegenerate.
\end{proof}

Next, we can find real numbers $a_i,b_j$ such that if $v = v_1 + a_1Jv_1 +
a_2J^2v_1 + a_3J^3v_1 + \cdots + a_{s-1}J^{s-1}v_1 + b_1Jw_1 + b_3J^3w_1 
+ b_5J^5w_1 + \cdots + b_{s-1}J^{s-1}w_1$, then
$$ \langle v,J^h v\rangle = \langle v,J^hw_1 \rangle = 0\qquad\mbox{for 
}\, 0\le h \le s-2\,.$$
Replacing $v_1$ by this $v$, we may assume that
\begin{eqnarray*}
\langle v_1, J^hv_1 \rangle &=& 0 \qquad\mbox{for }\, 0\le h\, ,\\
\langle v_1, J^hw_1 \rangle &=& 0 \qquad\mbox{for }\, 0\le h\le s-2\, ,\\
\langle J^{s-1}v_1, w_1 \rangle &=& 1\,.
\end{eqnarray*}
Now we can find real $c_i$ such that if $w = w_1 + c_1Jv_1 + c_3J^3v_1 +
c_5J^5v_1 + \cdots + c_{s-1}J^{s-1}v_1$, then $\langle w, J^hw \rangle =
0$ for all $h\ge 0$.  Replacing $w_1$ by this $w$, we may assume that all
inner products $\langle v_1, J^hv_1 \rangle$, $\langle v_1,J^hw_1
\rangle$, $\langle w_1,J^hw_1 \rangle$ vanish except $\langle
J^{s-1}v_1,w_1 \rangle =1$.  Setting $v_h = J^{h-1}v_1$ and $w_h =
J^{h-1}w_1$, the matrix of $J$ on the basis $\{v_1,\ldots ,v_s,w_1,\ldots
,w_s\}$ is $2s\times 2s$ of the form
$$ \arraycolsep .5em \left[ 
   \begin{array}{ccrccccrcc} 0 &&&&&&&&&\\ 
   1 & 0 &&&&&&&&\\[-1.1ex] &1& \ddots &&&&&&&\\ &&1& 0 &&&&&&\\ 
   &&& 1 & 0 &&&&&\\ &&&&& 0 &&&&\\ &&&&& 1 & 0 &&&\\[-1.1ex] 
   &&&&&&1& \ddots &&\\ &&&&&&&1& 0 &\\ &&&&&&&& 1 & 0\\
\end{array} \right] $$
and all inner products vanish except $\langle v_i,w_j \rangle = 
(-1)^{j-1}$ if $i+j=s+1$.

Lastly, assume that $s$ is odd. Consider $J^{s-1}v_1$ for $v_1\in V$. Then
there exists $w_1\in V$ such that $\langle J^{s-1}v_1, w_1\rangle = \half$
and $\langle J^{s-1}(v_1+w_1), v_1+w_1\rangle = 1$. We proceed by 
replacing $v_1$ by $v_1+w_1$. As before, we then have that $\lsp v_1, 
Jv_1, \ldots , J^{s-1}v_1 \rsp$ is a $J$-invariant, nondegenerate, 
$s$-dimensional subspace of $V$. Also, we can find real numbers $a_2, a_4,
\ldots , a_{s-1}$ such that 
$$ \langle v_1 + a_2J^2v_1 + \cdots + a_{s-1}J^{s-1}v_1, J^h \left( v_1 + 
a_2J^2v_1 + \cdots + a_{s-1}J^{s-1}v_1 \right) \rangle = 0 $$
for $0\le h\le s-2$. Replacing the current $v_1$ by $v_1 + a_2J^2v_1 + 
\cdots + a_{s-1}J^{s-1}v_1$, we obtain
\begin{eqnarray*}
\langle v_1, J^hv_1 \rangle &=& 0 \qquad\mbox{for }\, 0\le h\le s-2\, ,\\
\langle v_1, J^{s-1}v_1 \rangle &=& 1\,.
\end{eqnarray*}
Set $v_h = J^{h-1}v_1$ for $1\le h\le s$. The the matrix of $J$ with 
respect to the basis $\{v_1, v_2, \ldots , v_s\}$ is
$$ \arraycolsep .5em \left[
\begin{array}{ccrcc} 0 &&&&\\ 1 & 0 &&&\\[-1.1ex] &1& \ddots &&\\ 
   &&1& 0 &\\ &&& 1 & 0 \\ \end{array} \right] $$
with all inner products of basis vectors vanishing except $ \langle v_i, 
v_j \rangle = (-1)^{j-1}$ if $i+j=s+1$.

We finish almost as usual, with $1\le r\le s$, but $r\times r$ blocks now.
This completes the classification.

\end{document}